\documentclass[a4paper,12pt,reqno]{article}

\usepackage{amsthm}
\usepackage{amsmath,amssymb,latexsym,amsfonts,mathrsfs}
\usepackage[dvipdfmx]{graphicx}
\usepackage{bm}

\usepackage{color}
\usepackage{comment}
\usepackage{mathtools}
\usepackage{fancyhdr}
\usepackage[top=30truemm,bottom=30truemm,left=20truemm,right=20truemm]{geometry}
\usepackage{enumerate}

\newcommand{\R}{\mathbb{R}}
\newcommand{\C}{\mathbb{C}}
\newcommand{\N}{\mathbb{N}}
\newcommand{\Z}{\mathbb{Z}}

\newcommand{\SCR}[1]{{\mathscr #1}}

\newcommand{\CAL}[1]{{\cal #1}}

\newcommand{\J}[1]{\left\langle #1 \right\rangle}

\theoremstyle{plain}
\newtheorem{Thm}{{\bf Theorem}}[section]

\newtheorem{Prop}[Thm]{{\bf Proposition}}
\newtheorem{Cor}[Thm]{{\bf Corollary}}

\theoremstyle{definition}
\newtheorem{Def}[Thm]{{\bf Definition}}

\newcounter{Exami}

\allowdisplaybreaks[2]

\makeatletter
 \def\address#1#2{\begingroup
 \noindent\parbox[t]{7.8cm}{%
 \small{\scshape\ignorespaces#1}\par\vskip0ex
 \noindent\small{\itshape E-mail}%
 \/: #2\par\vskip4ex}\hfill%
 \endgroup}%
\makeatother

\begin{document}
\fontencoding{T1}\selectfont
%%%%%%%%%%%%%%%%%%%%%%%%%%%%%%%%%%%%%%%%%%%%%%%%%%%%%%%%%%%%%%% 
\begin{comment}
\end{comment}
%%%%%%%%%%%%%%%%%%%%%%%%%%%%%%%%%%%%%%%%%%%%%%%%%%%%%%%%%%%%%%%  
\title{Boundedness of propagators for Dirac equations with potentials on Wiener amalgam spaces}
\author{Shun Takizawa}
\date{\today}
\maketitle
\begin{abstract}
In this paper we prove boundedness of propagators for Dirac equations with unbounded time-dependent potentials on Wiener amalgam spaces. In particular we deal with class of potentials such as including Stark and harmonic potentials.
\end{abstract}
%%%%%%%%%%%%%%%%%%%%%%%%%%%%%%%%%%%%%%%%%%%%%%%%%%%%%%%%%%%%%%%%%%%%%%%%%%%%%%%%%%%%%
%%%%%%%%%%%%%%%%%%%%%%%%%%     セクション１（イントロ）      %%%%%%%%%%%%%%%%%%%%%%%%%%%%%%%%%%%%%
%%%%%%%%%%%%%%%%%%%%%%%%%%%%%%%%%%%%%%%%%%%%%%%%%%%%%%%%%%%%%%%%%%%%%%%%%%%%%%%%%%
\section{Introduction} 
In this paper, we study the following Cauchy problem for Dirac equations with potentials:
\begin{align}\label{CP1}
\begin{cases}
&i\partial_{t}u(t,x)= \CAL{D}u(t,x)+V(t,x)u(t,x), \hspace{3mm}(t,x)\in \R \times \R^d,
\\&u(0,x)=u_{0}(x),\hspace{3mm}x\in\R^d,
\end{cases}
\end{align}
where $u(t), u_0$ are $\C^n$-valued functions, the potential $V(t,x)$ is the square matrix of order $n$ with 
and $i=\sqrt{-1}$. Dirac operator $\CAL{D}$ is defined by
\begin{equation}\label{Diracop}
\CAL{D}=\alpha_0-i\sum_{j=0}^{d}\alpha_j \partial_{x_j}
\end{equation} 
where Hermitian matrices $\alpha_0, \ldots, \alpha_d \in \C^{n\times n}$ satisfy 
\begin{equation}\label{anticom}
\alpha_j \alpha_k+\alpha_k \alpha_j=2 \delta_{j, k} I_n.
\end{equation}
Here $\delta_{j, k}$ and $I_n$ stand for the Kronecker delta and identity matrix of order $n$, respectively.
For $d=3$ and $n=4$ the standard choice for $\alpha_0, \ldots, \alpha_d$ is the so-called Dirac representation:
\begin{align*}
\alpha_0=\begin{pmatrix}
   I_2 & 0 \\
   0 & -I_2
\end{pmatrix},\hspace{5mm} 
\alpha_j=\begin{pmatrix}
   0 & a_j \\
   a_j & 0
\end{pmatrix},
\hspace{3mm}j=1,2,3, 
\end{align*}
where $a_j$ are Pauli matrices:
\begin{align*}
a_1=\begin{pmatrix}
   0 & 1 \\
   1 & 0
\end{pmatrix}, \hspace{3mm}
a_2=\begin{pmatrix}
   1 & 0 \\
   0 & -1
\end{pmatrix}, \hspace{3mm}
a_3=\begin{pmatrix}
 0 & -i \\
   i & 0
\end{pmatrix}.
\end{align*}
In general, for any $d$ there exist $n=n(d)$ and $\alpha_0,\ldots, \alpha_d \in \C^{n\times n}$ satisfying (\ref{anticom}) (see Kalf and Yamada \cite{Kalf-Yamada}).

There are numerous related works.
First we remark on known results in cases where potentials are absent.
H\"{o}rmander \cite{Hormander} proved that $e^{i|D|^2}$ is bounded on $L^p$ if and only if $p=2$.
Otherwhile on modulation spaces $M^{p,q}$ (see Definition \ref{DefM} below), the operator $e^{i|D|^2}$ is bounded for all indices $1\leq p, q\leq \infty$ (see Gr\"{o}chenig and Heil \cite{Grochenig-Heil}, Toft \cite{Toft}, Wang, Zhao and Guo \cite{WZG}). More generally B\'{e}nyi, Gr\"{o}chenig Okoudjou and Rogers \cite{BGOR} studied the boundedness of $e^{i|D|^{\kappa}}$ for $0\leq\kappa\leq2$ on $M^{p,q}$ for all $1\leq p, q\leq\infty$.
Cunanan and Sugimoto \cite{Cunanan-Sugimoto} demonstrated that boundedness of $e^{i|D|^{\kappa}}$  for $0\leq\kappa\leq1$ on Wiener amalgam spaces $W^{p,q}$ (see Definition \ref{DefWA} below) for all $1\leq p, q\leq\infty$. 
Boundedness  of $e^{i|D|^{\kappa}}$ for $\kappa\geq 1$ from $W_{0,\rho}^{p,q}$ to $W^{p,q}$ has studied by T. Kato and Tomita \cite{Kato-Tomita} and Guo and Zhao \cite{Guo-Zhao}.

Next we remark on previous works in cases where potentials are present.
On Schr\"{o}dinger equations with potentials, if a potential is at most linear growth, then the corresponding propagator is bounded on $M^{p,q}$ for all $1\leq p, q\leq \infty$ and if a potential is at most quadratic growth, then the corresponding propagator is bounded on $M^{p,p}$ (i.e. $p=q$) (see Cordero and Nicola \cite{CN}, K. Kato, Kobayashi and Ito \cite{KKI}, Cordero, Nicola and Rodino \cite{CNR}).
On replacing Schr\"{o}dinger equations with Dirac equations, similar results hold (see K. Kato and Naumkin \cite{Kato-Naumkin}, Trapasso \cite{Trapasso1}). As far as we know, there are no papers concerning boundedness of propagators for Dirac equation with unbounded potentials in the framework of Wiener amalgam spaces. 
%
%%%%%%%%%%%%%%%%%%%%%%%%       主定理      %%%%%%%%%%%%%%%%%%%%%%%%%%%%
%
%

The following theorems are our main results, which argue boundedness of propagators with at most linear and quadratic growth potentials.
\begin{Thm}\label{thm1}
Let $1\leq p, q\leq \infty, \rho\in\R$ and $T>0$. Suppose that $V$ satisfies $V(t,x)\equiv V(x)$, 
\begin{equation}
\partial_{x}^{\alpha}V\in M_{0,|\rho|}^{\infty,1}(\R^d; \R^{n\times n}) \hspace{3mm}\text{with} \hspace{3mm}|\alpha|=1
\end{equation}
or $V(t,x)\equiv Q(t,x)I_n$,
\begin{equation}
\partial_{x}^{\alpha}Q\in C_b(\R; M_{0,|\rho|}^{\infty,1}(\R^d)) \hspace{3mm}\text{with} \hspace{3mm}|\alpha|=1.
\end{equation}
For every $u_0 \in \CAL{W}_{r,\rho}^{p,q}(\R^d;\C^n))$, there exists a unique solution $u\in C([0,T]; \CAL{W}_{r,\rho}^{p,q}(\R^d;\C^n))$ to (\ref{CP1}). Moreover the corresponding propagator $U(t)$ is uniformly bounded on $\CAL{W}_{r,\rho}^{p,q}(\R^d;\C^n))$ with respect to $t\in (-T, T)$.
\end{Thm}
\begin{Thm}\label{thm2}
Let $1\leq p, q\leq \infty, \rho\in\R$ and $T>0$. Suppose that $V$ satisfies $V(t,x)=Q(t,x)I_n$,
\begin{equation}
\partial_{x}^{\alpha}Q\in C_b(\R; M^{\infty,1}(\R^d)) \hspace{3mm}\text{with} \hspace{3mm}|\alpha|=2.
\end{equation}
For every $u_0 \in \CAL{W}_{r,0}^{p,q}(\R^d;\C^n)$, there exists a unique solution $u\in C([0,T]; \CAL{W}_{r,0}^{p,q}(\R^d;\C^n))$ to (\ref{CP1}). Moreover the corresponding propagator $U(t)$ is uniformly bounded on $\CAL{W}_{r,0}^{p,q}(\R^d;\C^n))$ with respect to $t\in (-T, T)$.
\end{Thm}

As a related result of Theorem \ref{thm1} and \ref{thm2} in the present paper, there is Theorem 1.2 in \cite{Trapasso1}, which is shown that roughly speaking, boundedness of propagators of Dirac equations with bounded potentials on $W^{p,q}$ for all $1\leq p, q\leq \infty$.
The contribution of the present paper is to give boundedness of the propagators on $W^{p,q}$ with unbounded potentials for all indices $1\leq p, q\leq\infty$. Typical examples of assumptions in  Theorem \ref{thm1} and Theorem \ref{thm2} are Stark and harmonic potentials, respectively.
%%%%%%%%%%%%%%%%%%%%%%%%%%%%%%%%%%%%%%%%%%%%%%%%%%%%%%%%%%%%%%%%%%%%%%%%%%%%%%%%%%%%
%%%%%%%%%%%%%%%%%%%%%%%%%%　セクション2 (予備概念）     %%%%%%%%%%%%%%%%%%%%%%%%%%%%%%%%%%%
%%%%%%%%%%%%%%%%%%%%%%%%%%%%%%%%%%%%%%%%%%%%%%%%%%%%%%%%%%%%%%%%%%%%%%%%%%%%%%%%%%%%
\section{Preliminaries}
\subsection{Notation}
Let $\nabla_x=(\partial_{x_1},\ldots,\partial_{x_d})$ and $\J{\cdot}=\left( 1+|\cdot|^2\right)^{1/2}$.  We set $|x|_{\C^m}=\max\{|x_1|, \ldots, |x_m|\}$ for $x=(x_1,\ldots, x_m)\in \C^m$.
We use notation $\SCR{F}f(\xi)=\widehat{f}(\xi)=(2\pi)^{-d} \int_{\R^d}f(x)e^{-ix\cdot\xi}dx$ for Fourier transform of $f$ 
and $\SCR{F}^{-1}f(x)=\check{f}(x)=\int_{\R^d}f(\xi)e^{ix\cdot\xi}d\xi$ for the  inverse Fourier transform of $f$. We often write $\int$ instead of $\int_{\R^m}$ with $m\in \N$ for short.
We write $X\lesssim Y$ if $X\leq CY$ with some constant $C>0$ in the proofs.
We define
\begin{equation}
C_{\geq k}^{\infty}(\R^d; \C^n)=\{f\in C^{\infty}(\R^d; \C^n): |\partial^{\alpha}f|_{\C^n}\leq C_{\alpha} \hspace{2mm}\text{with}\hspace{2mm} |\alpha|\geq k\}
\end{equation}
and $C_b(\R; X)$ denotes the space of continuous and bounded function $\R \to X$.
\subsection{Short-time Fourier transform (STFT)} 
We introduce the vector-valued Short-time Fourier transform used difinitions of modulation and Wiener amalgam spaces.
\begin{Def}[STFT]
Let $g \in \CAL{S}(\mathbb{R}^d)\setminus \{0\}$ and $f \in \CAL{S} '(\R^d; \C^n)$. Then the short-time Fourier transform $V_{g}f$ of $f$ with respect to a window function $g$ is defined by
\begin{equation*}
V_{g}f(x, \xi)=\int_{{\mathbb R}^d} \overline{g(y-x)}f(y)e^{-iy\cdot \xi}dy,\quad (x, \xi) \in{\mathbb R}^d\times{\mathbb R}^d.
\end{equation*}
We also define the formal adjoint operator $V_{g}^{*}$ of $V_{g}$ by
\begin{equation*}
V_{g}^{*}F(x)=\iint_{{\mathbb R}^{2d}} g(x-y)F(y, \xi)e^{ix\cdot \xi}dy\bar{d}\xi,\quad x\in{\mathbb R}^d, F\in \CAL{S}^{'}(\R^{2d}; \C^n)
\end{equation*}
with $\bar{d}\xi=(2\pi)^{-d}d\xi$. 
\end{Def}
\begin{Prop}[Inverse formula for STFT]\label{Prop1}
Let $g \in \CAL{S}(\mathbb{R}^d)\setminus \{0\}$. Then it follows that
\begin{equation*}
f=\|g\|_{L^2}^{-2} V_{g}^{*}[V_{g}f], \hspace{3mm} f\in \CAL{S}'(\R^d; \C^n).
\end{equation*}
\end{Prop}
Proposition \ref{Prop1} is proved in \cite[Proposition 2.5]{Wahlberg}.
\subsection{Modulation and Wiener amalgam spaces}
\begin{Def}[Modulation spaces]\label{DefM}
Let $1\leq p, q\leq \infty$ and $g\in \CAL{S}(\R^n)\setminus \{0\}$. Then we define
\begin{align*}
&M_{r,\rho}^{p,q}(\R^d; \C^n)=\left\{f\in\CAL{S}'(\R^d) : \left\| \|\J{x}^r \J{\xi}^{\rho}  |V_{g}f(x,\xi)|_{\C^n}\|_{L_{x}^{p}}\right\|_{L_{\xi}^{q}}<\infty\right\},
\\
&\|f\|_{M_{r,\rho}^{p,q}(\R^d; \C^n)}=\left\| \| \J{x}^r \J{\xi}^{\rho} |V_{g}f(x,\xi)|_{\C^n}\|_{L_{x}^{p}}\right\|_{L_{\xi}^{q}}.
\end{align*}
\end{Def}
\begin{Def}[Wiener amalgam spaces]\label{DefWA}
Let $1\leq p, q\leq \infty$ and $g\in \CAL{S}(\R^n)\setminus \{0\}$. Then we define
\begin{align*}
&W_{r,\rho}^{p,q}(\R^d; \C^n)=\left\{f\in\CAL{S}'(\R^d) : \left\| \| \J{x}^r \J{\xi}^{\rho} |V_{g}f(x,\xi)|_{\C^n}\|_{L_{\xi}^{q}}\right\|_{L_{x}^{p}}<\infty\right\},
\\
&\|f\|_{W_{r,\rho}^{p,q}(\R^d; \C^n)}=\left\|  \| \J{x}^r \J{\xi}^{\rho} |V_{g}f(x,\xi)|_{\C^n}\|_{L_{\xi}^{q}}\right\|_{L_{x}^{p}}.
\end{align*}
\end{Def}
Modulation and Wiener amalgam spaces are independent of the choice of the window function $g$.
We write $M^{p,q}(\R^d; \C^n)$ and $W^{p,q}(\R^d; \C^n)$ instead of $M_{0,0}^{p,q}(\R^d; \C^n)$ and $W_{0,0}^{p,q}(\R^d; \C^n)$, respectively. We often write $M_{r,\rho}^{p,q}=M_{r,\rho}^{p,q}(\R^d; \C^n)$ and $W_{r,\rho}^{p,q}=W_{r,\rho}^{p,q}(\R^d; \C^n)$ for short.
We collect properties of modulation and Wiener amalgam spaces in order to prove Theorem \ref{thm1} and \ref{thm2}.
\begin{Prop}[\cite{Trapasso1}, Remark 2.6]
Let 
\begin{equation*}
\CAL{M}_{r,\rho}^{p,q}:=\overline{\CAL{S}}^{\|\cdot\|_{M_{r,\rho}^{p,q}}},
\hspace{3mm} \CAL{W}_{r,\rho}^{p,q}:=\overline{\CAL{S}}^{\|\cdot\|_{W_{r,\rho}^{p,q}}}.
\end{equation*}
Then it holds that 
\begin{equation*}
\CAL{M}_{r,\rho}^{p,q}=M_{r,\rho}^{p,q},
\hspace{3mm} \CAL{W}_{r,\rho}^{p,q}=W_{r,\rho}^{p,q},
\end{equation*}
for $1\leq p, q <\infty$.
\end{Prop}

\begin{Prop}\label{PropS}
Let $(D_{\theta}f)(x)=f(\theta x)$ for $0<\theta\leq1$ and $(T_{a}f)(x)=f(x-a)$ for $a\in \R^d$.
If $f\in M^{\infty,1}(\R^d; \C^n)$, then $D_{\theta}f,  T_{a}f\in M^{\infty,1}(\R^d; \C^n)$. 
In particular it holds that there exists $C>0$ such that for all $\theta\in(0,1]$ and $a\in \R^d$,
\begin{align}
&\|D_{\theta}f\|_{M^{\infty,1}}\leq C \|f\|_{M^{\infty,1}}, \label{Dtheta}
\\
&\|T_{a}f\|_{M^{\infty,1}}=\|f\|_{M^{\infty,1}}. \label{Ta}
\end{align}
\begin{proof}
The inequality (\ref{Dtheta}) follows from Sugimoto and Tomita \cite[Theorem 1.1]{Sugimoto-Tomita}.
It holds (\ref{Ta}) by Definition \ref{DefM}.
\end{proof}
%%%%%%%%%%%%%%%%
\begin{Cor}
Let $f\in M^{\infty,1}(\R^d; \C^n)$ and 
\begin{equation*}
\widetilde{f}(x,y):=\int_{0}^{1}f(x+\theta(y-x))(1-\theta)d\theta.
\end{equation*}
Then there exists constant $C>0$ independent of $x\in \R^d$ such that
\begin{equation*}
\|\widetilde{f}(x,\cdot)\|_{M^{\infty,1}}\leq C \|f\|_{M^{\infty,1}}.
\end{equation*}
\end{Cor}
%%%%%%%%%%%%%%%%
\end{Prop}
The following statement is complex interpolation theorem of Wiener amalgam spaces.
\begin{Prop}[\cite{Cunanan-Sugimoto}, Lemma 2.1]\label{Prop4}
Let $r,\rho\in\R$, $0<\theta<1$ and $p,q,p_j,q_j\in [1,\infty]$ for $j=1,2$. If $1/p=\theta/p_1+(1-\theta)/p_2$ and $1/q=\theta/q_1+(1-\theta)/q_2$, then
\begin{equation*}
\left(W_{r,\rho}^{p_1,q_1}, W_{r,\rho}^{p_2,q_2}\right)_{[\theta]}=W_{r,\rho}^{p,q}.
\end{equation*}
\end{Prop}
\begin{Prop}\label{prop210}
Let $\rho\in \R, k\in \N$. If $f: \R\times \R^d \to \C^n$ is $\partial^{\alpha}f\in C_b(\R; M_{0,|\rho|}^{\infty,1}(\R^d; \C^n))$ for any $\alpha\in \Z_{+}^{d}$ with $|\alpha|=k$, then there exist $f_1\in  C_b (\R; C_{\geq k}^{\infty}(\R^d; \C^n))$ and $f_2\in C_b (\R; M_{0,|\rho|}^{\infty,1}(\R^d; \C^n))$ such that $f=f_1+f_2$.
\end{Prop}
The proof of Proposition \ref{prop210} is similar to one of \cite[ Proposition 3.2]{Trapasso1}. 
%%%%%%%%%%%%%%%%%%%%%%%%%%%%%%%%%%%%%%%%%%%%%%%%%%%%%%%%%%%%%%%%%%%%%%%%%%%%%%%%%%
%%%%%%%%%%%%%%%%%%%%%%%%%%%%%%%%%%%%%%%%%%%%%%%%%%%%%%%%%%%%%%%%%%%%%%%%%%%%%%%%%%
%%%%%%%%%%%%%%%%%%%%%%%%%%%%%%%%%%%%%%%%%%%%%%%%%%%%%%%%%%%%%%%%%%%%%%%%%%%%%%%%%%
\section{Proof of Theorem \ref{thm1}}
\subsection{Construction of a parametrix}
First we consider a more general Fourier multiplier $\sigma(D)$ instead of Dirac operator $\CAL{D}$.
Suppose that the symbol $\sigma(\xi)$ of $\sigma(D)=\SCR{F}^{-1} \sigma(\xi) \SCR{F}$ is an Hermitian ($n\times n$)-matrix valued  and satisfies 
\begin{equation}\label{symbol}
|\partial_{\eta}^{\alpha} \left(\sigma(\xi+\eta)-\sigma(\xi)\right)|_{\C^{n\times n}}\leq C_{\alpha}\J{\eta}^k, \hspace{3mm}\xi,\eta\in \R^d,
\end{equation}
for some $k\geq 0$. The Dirac operator (\ref{Diracop}) satisfies (\ref{symbol}).
Moreover we consider the following Cauchy problem:
\begin{align}\label{CP2}
\begin{cases}
&i\partial_{t}u(t,x)= \sigma(D)u(t,x)+V(t,x)u(t,x), \hspace{3mm}(t,x)\in \R \times \R^d,
\\&u(0,x)=u_{0}(x),\hspace{3mm}x\in\R^d,
\end{cases}
\end{align}
instead of (\ref{CP1}).
By Proposition \ref{prop210},  we can decompose $V=V_1+V_2$ with $V_1\in C_b(\R; C_{\geq1}^{\infty}(\R^d; \C^{n\times n}))$ and $V_2 \in C_b(\R; M_{0,|\rho|}^{\infty,1}(\R^d; \C^{n\times n}))$.

We employ an approximate propagator which is similar to one used in Tataru \cite{Tataru}, Cordero, Nicola and Rodino \cite{CNR}.
We set
\begin{align*}
U_1(t,s)f(y)&=\iint e^{-i\int_{s}^{t}\left( \sigma(\xi)+V_1(\tau,x)\right)d\tau}e^{i(y-x)\xi}
 g(y-x)V_{g}f(x,\xi)dxd\xi
\end{align*}
for $f\in \CAL{S}(\R^d; \C^n)$ and $g\in \CAL{S}(\R^d)\setminus \{0\}$ with $\|g\|_{L^2}=1$. 
By Proposition \ref{Prop1}, it holds
\begin{equation}\label{unitary}
U_1(s,s)f=f, \hspace{3mm} f\in \CAL{S}(\R^d; \C^n).
\end{equation}
By $V_{g}f\in \CAL{S}(\R^{2d})$ (see Be\'{n}yi and Okoudjou \cite{BO}), we can apply dominated convergence theorem and obtain
\begin{equation}\label{id13}
i\partial_{t}U_1(t,s)f(y)
=\iint ( \sigma(\xi)+V_1(t,x)) e^{-i\int_{s}^{t}\left( \sigma(\xi)+V(\tau,x)\right)d\tau}e^{i(y-x)\xi} V_{g}f(x,\xi)dx d\xi.
\end{equation}
We obtain
\begin{align*}
&\sigma(D_y)[g(\cdot-x)e^{i(\cdot-x)\xi}](y)
\\
&=\int e^{iy\cdot\eta}\sigma(\eta)\SCR{F}[g(\cdot-x)e^{i(\cdot-x)\xi}](\eta)d\eta
\\
&=\int e^{iy\cdot\eta}\sigma(\eta) e^{-ix\eta} \widehat{g}(\eta-\xi)d\eta
\\
&=e^{i(y-x)\xi} \int e^{i(y-x)\eta} \sigma(\eta+\xi) \widehat{g}(\eta)d\eta,
\end{align*}
which implies
\begin{align}\label{id14}
&\sigma(D_y)U_0(t,s)f(y) \notag
\\
&=\iint e^{-i\int_{s}^{t}\left( \sigma(\xi)+V_1(\tau,x)\right)d\tau} e^{i(y-x)\xi} \left(\int e^{i(y-x)\eta} \sigma(\eta+\xi) \widehat{g}(\eta)d\eta\right)
 V_{g}f(x,\xi)dxd\xi.
\end{align} 
By the Taylor expansion, the identity
\begin{equation}\label{id15}
V(t,y)=V_1(t,x)+\sum_{|\alpha|=1} \int_{0}^{1}\left(\partial_{x}^{\alpha}V_1\right)(t,x+\theta(y-x)) (1-\theta)d\theta (y-x)^{\alpha}+V_2(t,y)
\end{equation}
holds. From (\ref{id13}), (\ref{id14}) and (\ref{id15}), we have
\begin{align*}
&\left(i\partial_{t}-\sigma(D_{y})-V(t,y)\right)U_1(t,s)f(y)
\\
&=-\iint e^{-i\int_{s}^{t}\left( \sigma(\xi)+V_1(\tau,x)\right)d\tau}e^{i(y-x)\xi} V_{g}f(x,\xi)(R_{\sigma}(t,x,y,\xi)+R_{V1}(t,x,y)+R_{V_2}(t,x,y)) dx d\xi, 
\end{align*}
where 
\begin{align*}
&R_{\sigma}(t,x,y,\xi) =\int e^{i(y-x)\eta} \left( \sigma(\eta+\xi)-\sigma(\xi) \right) \widehat{g} (\eta) d\eta,
\\
&R_{V_1}(t,x,y)=\sum_{|\alpha|=1} \int_{0}^{1}\left(\partial_{x}^{\alpha}V_1\right)(t,x+\theta(y-x)) (1-\theta)d\theta (y-x)^{\alpha} g(y-x),
\\
&R_{V_2}(t,x,y)=V_{2}(t,y)g(y-x).
\end{align*}
Let
\begin{align*}
&K_{\sigma}(t,s)f(y)=-\iint e^{-i\int_{s}^{t}\left( \sigma(\xi)+V_1(\tau,x)\right)d\tau}e^{i(y-x)\xi} V_{g}f(x,\xi) R_{\sigma}(t,x,y,\xi) dx d\xi,
\\
&K_{V_j}(t,s)f(y)= -\iint e^{-i\int_{s}^{t}\left( \sigma(\xi)+V_1(\tau,x)\right)d\tau}e^{i(y-x)\xi} V_{g}f(x,\xi) R_{V_j}(t,x,y) dx d\xi, \hspace{3mm}j=1,2.
\end{align*}
Then we have
\begin{equation}\label{DE}
\left(i\partial_{t}-\sigma(D_{y})-V(t,y)\right)U_1(t,s)f(y)=K_{\sigma}(t,s)f(y)+K_{V_1}(t,s)f(y)+K_{V_2}(t,s)f(y).
\end{equation}
%%%%%%%%%%%%%%%%%%%%%%%%%%%%%%%%%%%%%%%%%%%%%%%%%%%%%%%%%%
%%%%%%%%%%%%%%%%%%%%%%%%%%%%%%%%%%%%%%%%%%%%%%%%%%%%%%%%%%
\subsection{Proof of boundedness of operators}
We prove that operators $U_1(t,s), K_{\sigma}(t,s), K_{V_1}(t,s)$ and $K_{V_2}(t,s)$ are uniformly bounded on $W_{r,\rho}^{p,q}$ with respect to $0\leq|s-t|\leq T$ for all $1\leq p, q\leq\infty$ and $r, \rho \in \R$.
First we prove the following inequality:
\begin{equation} \label{inqU1}
\|U_1(t,s)f\|_{W_{r,\rho}^{p,q}}\lesssim \|f\|_{W_{r,\rho}^{p,q}}.
\end{equation}
Let $\varphi \in \CAL{S}(\R^d)\setminus \{0\}$ be fixed.
Using the equality $(1-\Delta_{y})^{N}e^{iy(\xi-\Xi)}=\J{\xi-\Xi}^{2N}e^{iy(\xi-\Xi)}$, we have for any $N\in \N,$
\begin{align} \label{inq09}
&|V_{\varphi}[U_1(t,s)f](X,\Xi)|_{\C^n} \notag
\\
&=\left|\int_{\R^{3d}}\overline{\varphi(y-X)}e^{-i\int_{s}^{t}\left( \sigma(\xi)+V_1(\tau,x)\right)d\tau}e^{i(y-x)\xi-iy\Xi} g(y-x) V_{g}f(x,\xi)dxdyd\xi \right| \notag
\\
&\lesssim \sum_{|\beta_1+\beta_2|\leq 2N} \int_{\R^{3d}} \J{\xi-\Xi}^{-2N} |\partial_{y}^{\beta_1}\varphi(y-X)| |\partial_{y}^{\beta_2}g(y-x)| |V_{g}f(x,\xi)|_{\C^n} dx dy d\xi.
\end{align}
By the Peetre inequalities $\J{X}^{r}\lesssim \J{X-y}^{|r|} \J{y-x}^{|r|} \J{x}^{r}$ and $\J{\Xi}^{\rho}\lesssim \J{\Xi-\xi}^{|\rho|} \J{\xi}^{\rho}$, we obtain 
\begin{align} \label{inq11}
&\|\|\J{X}^{r} \J{\Xi}^{\rho} |V_{\varphi}[U_1(t,s)f](X,\Xi)|_{\C^n} \|_{L_{\Xi}^{q}} \|_{L_{X}^{p}} \notag
\\
&\lesssim \sum_{|\beta_1+\beta_2|\leq 2N} \left\| \left\| \int_{\R^{3d}} \J{\xi-\Xi}^{-2N+|\rho|} |\widetilde{\varphi}(y-X) \widetilde{g}(y-x)|  \J{x}^{r} \J{\xi}^{\rho} |V_{g}f(x,\xi)|_{\C^n} dx dy d\xi \right\|_{L_{\Xi}^{q}} \right\|_{L_{X}^{p}},
\end{align}
where $\widetilde{\varphi}(x)=\J{x}^{|r|} \partial_{x}^{\beta_1}\varphi(x)$ and $\widetilde{g}(x)=\J{x}^{|r|} \partial_{x}^{\beta_2}g(x)$.
In particular for $N\in \N$ with $2N>n+|\rho|$ we find by Minkowski's inequality and Fubini's theorem,
\begin{align} \label{inq12}
&\left\| \left\| \int_{\R^{3d}} \J{\xi-\Xi}^{-2N+|\rho|} |\widetilde{\varphi}(y-X) \widetilde{g}(y-x)|  \J{x}^{r} \J{\xi}^{\rho} |V_{g}f(x,\xi)|_{\C^n} dx dy d\xi \right\|_{L_{\Xi}^{1}} \right\|_{L_{X}^{1}} \notag
\\
&\leq  \int_{\R^{3d}}  \left\|\widetilde{\varphi}(y-X) \right\|_{L_{X}^{1}} | \widetilde{g}(y-x)| \J{x}^{r} \J{\xi}^{\rho} |V_{g}f(x,\xi)|_{\C^n}  \left\| \J{\xi-\Xi}^{-2N+|\rho|} \right\|_{L_{\Xi}^{1}}  dx dy d\xi \notag
\\
&\lesssim  \int_{\R^{2d}}  \left( \int |\widetilde{g}(y-x)| dy \right) \J{x}^{r} \J{\xi}^{\rho} |V_{g}f(x,\xi)|_{\C^n}   dx d\xi \notag
\\
&\lesssim \|\|\J{x}^{r} \J{\xi}^{\rho} |V_{g}f(x,\xi)|_{\C^n} \|_{L_{\xi}^{1}} \|_{L_{x}^{1}}. 
\end{align}
(\ref{inq11}) together with (\ref{inq12}) implies
\begin{equation} \label{inq3}
\|U_1(t,s)f\|_{W_{r,\rho}^{1,1}} \lesssim \|f\|_{W_{r,\rho}^{1,1}}.
\end{equation} 
Similary we have
\begin{equation} \label{inq4}
\|U_1(t,s)f\|_{W_{r,\rho}^{p,q}} \lesssim \|f\|_{W_{r,\rho}^{p,q}}
\end{equation}
with $(p,q)=(1,\infty), (\infty,1), (\infty, \infty)$.
By (\ref{inq3}), (\ref{inq4}) and Proposition \ref{Prop4}, (\ref{inqU1}) holds.

Next we show 
\begin{equation} \label{Ksigma}
\|K_{\sigma}(t,s)f\|_{W_{r,\rho}^{p,q}}\lesssim \|f\|_{W_{r,\rho}^{p,q}}.
\end{equation}
We note that 
\begin{equation} \label{Rsigma} 
|\partial_{y}^{\alpha}R_{\sigma}(t,x,y,\xi)|\leq C_{\alpha, N} \J{y-x}^{-2N}
\end{equation}
for all $N \in \N, \alpha \in \Z_{+}^{n}$.
Indeed, it follows from $g\in \CAL{S}$ and (\ref{symbol}) that 
\begin{align*}
&|\partial_{y}^{\alpha}R_{\sigma}(t,x,y,\xi)|=\left|\int (i\eta)^{\alpha} e^{i(y-x)\eta}  \left(\sigma(\eta+\xi)-\sigma(\xi) \right) \widehat{g} (\eta) d\eta\right|
\\
&=\J{y-x}^{-2N} \left|\int e^{i(y-x)\eta} (1-\Delta_{\eta})^{N}  \left[\left(\sigma(\eta+\xi)-\sigma(\xi) \right) \eta^{\alpha} \widehat{g} (\eta) \right] d\eta\right|
\\
&\lesssim \J{y-x}^{-2N}
\end{align*}
(see also \cite{Kato-Naumkin}).
Since (\ref{Rsigma}) holds, we can also estimate replacing $g(y-x)$ with $R_{\sigma}(t,x,y,\xi)$ in (\ref{inq09})-(\ref{inq4}). Thus we derive (\ref{Ksigma}) by a similar argument to (\ref{inqU1}).

Thirdly we demonstrate 
\begin{equation} \label{Kv1}
\|K_{V_1}(t,s)f\|_{W_{r,\rho}^{p,q}}\lesssim \|f\|_{W_{r,\rho}^{p,q}}.
\end{equation}
We write $V_{1}^{\alpha}(t,x,y)=\int_{0}^{1}\left(\partial_{x}^{\alpha}V_1\right)(t,x+\theta(y-x)) (1-\theta)d\theta$ and $g^{\alpha}(x)=x^{\alpha} g(x)$ for short.
Then $R_{V1}(t,x,y)=\sum_{|\alpha|=1}V_{1}^{\alpha}(t,x,y) g^{\alpha}(y-x)$ is valid.
By the equality $\J{\xi-\Xi}^{2N}e^{iy(\xi-\Xi)}=(1-\Delta_y)^{N}e^{iy(\xi-\Xi)}$ for $N\in \N$ and $V_1\in C_b(\R; C_{\geq1}^{\infty}(\R^d; \C^{n\times n}))$, we have
\begin{align} \label{inq100}
&|V_{\varphi}[K_{V_1}(t,s)f](X,\Xi)|_{\C^n} \notag
\\
&\leq \sum_{|\alpha|=1} \left| \int_{\R^{3d}} \overline{\varphi(y-X)} e^{-i\int_{s}^{t}\left( \sigma(\xi)+V_1(\tau,x)\right)d\tau}e^{i(y-x)\xi-iy\Xi} g^{\alpha}(y-x) V_{1}
^{\alpha}(t,x,y) V_{g}f(x,\xi) dxdyd\xi \right| \notag
\\
&\lesssim \sum_{|\alpha|=1}\sum_{|\beta_1+\beta_2+\beta_3|\leq 2N}\int_{\R^{3d}} \J{\xi-\Xi}^{-2N} |\partial^{\beta_1}\varphi(y-X) \partial^{\beta_2}g^{\alpha}(y-x)|  \notag
\\
& \hspace{50mm} \times |(\partial_{y}^{\beta_3} V_{1}^{\alpha})(t,x,y)|_{\C^{n\times n}} |V_{g}f(x,\xi)|_{\C^n} dxdyd\xi \notag
\\
&\lesssim \sum_{|\alpha|=1}\sum_{|\beta_1+\beta_2|\leq 2N}\int_{\R^{3d}} \J{\xi-\Xi}^{-2N} |\partial^{\beta_1}\varphi(y-X) \partial^{\beta_2}g^{\alpha}(y-x)| |V_{g}f(x,\xi)|_{\C^n} dxdyd\xi.
\end{align}
Taking $\|\|\J{X}^r \J{\Xi}^{\rho} \cdot\hspace{1mm}\|_{L_{\Xi}^q}\|_{L_{X}^p}$ both sides of (\ref{inq100}), we obtain
(\ref{Kv1}) by the same argument as (\ref{inq4}).

Finally we show  
\begin{equation} \label{Kv2}
\|K_{V_2}(t,s)f\|_{W_{r,\rho}^{p,q}}\lesssim \|f\|_{W_{r,\rho}^{p,q}}.
\end{equation}
Inverse formula for STFT and integration by parts follow that 
\begin{align} \label{inq5}
&|V_{\varphi}[K_2(t,s)f](X,\Xi)| \notag
\\
&=\left| \int_{\R^{5d}} \overline{\varphi(y-X)} g(y-x) \psi(y-z) e^{-i\int_{s}^{t}\left( \sigma(\xi)+V_1(\tau,x)\right)d\tau-ix\xi}e^{iy(\xi+\eta-\Xi)} \right. \notag
\\
&\left. \hspace{25mm} \times V_{\psi}[V_2(t,\cdot)](z,\eta) V_{g}f(x,\xi)dxdydzd\xi d\eta \right| \notag
\\
&\lesssim \sum_{|\beta_1+\beta_2+\beta_3|\leq 2N} \int_{\R^{5d}} \J{\xi+\eta-\Xi}^{-2N} \left| \partial_{y}^{\beta_1}\overline{\varphi(y-X)} \partial_{y}^{\beta_2}g(y-x) \partial_{y}^{\beta_3}\psi(y-z) \right| \notag
\\
& \hspace{40mm} \times \left|V_{\psi}[V_2(t,\cdot)](z,\eta) V_{g}f(x,\xi) \right| dxdydzd\xi d\eta
\end{align}
for $\psi\in \CAL{S}(\R^n)$ with $\|\psi\|_{L^2}=1$ and $N\in \N$.
Combining the Peetre inequalities $\J{X}^{r}\lesssim \J{X-y}^{|r|} \J{y-x}^{|r|} \J{x}^{r}$, $\J{\Xi}^{\rho}\lesssim   \J{\Xi-\xi-\eta}^{|\rho|}\J{\eta}^{|\rho|} \J{\xi}^{\rho}$ and (\ref{inq5}), we observe
\begin{align} \label{inq6}
&\J{X}^r \J{\Xi}^{\rho} |V_{\varphi}[K_2(t,s)f](X,\Xi)| \notag
\\
&\lesssim \sum_{|\beta_1+\beta_2+\beta_3|\leq 2N} \int_{\R^{5d}} \J{\xi+\eta-\Xi}^{-2N+|\rho|} \left| \widetilde{\varphi}(y-X) \widetilde{g}(y-x) \widetilde{\psi}(y-z) \right| \notag
\\
& \hspace{40mm} \times \J{\eta}^{|\rho|} |V_{\psi}[V_2(t,\cdot)](z,\eta)| \J{x}^r \J{\xi}^{\rho}|V_{g}f(x,\xi) | dxdydzd\xi d\eta,
\end{align}
where $\widetilde{\varphi}(x)=\J{x}^{|r|}\partial_{x}^{\beta_1}\varphi(x)$, $\widetilde{g}(x)=\J{x}^{|r|}\partial_{x}^{\beta_2}g(x)$ and $\widetilde{\psi}(x)=\partial_{x}^{\beta_3}\psi(x)$.
We have 
\begin{align} \label{inq7}
&\left\| \left\|\int_{\R^{5d}} \J{\xi+\eta-\Xi}^{-2N+|\rho|} \left| \widetilde{\varphi}(y-X) \widetilde{g}(y-x) \widetilde{\psi}(y-z) \right| \right.\right. \notag
\\
& \hspace{20mm} \times \left.\left. \J{\eta}^{|\rho|} |V_{\psi}[V_2(t,\cdot)](z,\eta)| \J{x}^r \J{\xi}^{\rho}|V_{g}f(x,\xi) | dxdydzd\xi d\eta \right\|_{L_{\Xi}^{1}} \right\|_{L_{X}^{1}} \notag
\\
&\lesssim \int_{\R^{5d}} \left\| \J{\xi+\eta-\Xi}^{-2N+|\rho|} \right\|_{L_{\Xi}^{1}} \left\|\widetilde{\varphi}(y-X) \right\|_{L_{X}^{1}}  |\widetilde{g}(y-x) | \J{x}^r \J{\xi}^{\rho}|V_{g}f(x,\xi) |   \notag
\\
& \hspace{20mm} \times  \J{\eta}^{|\rho|} |V_{\psi}[V_2(t,\cdot)](z,\eta) \widetilde{\psi}(y-z)| dxdydzd\xi d\eta  \notag
\\
&\lesssim \int_{\R^{4d}}   |\widetilde{g}(y-x) | \J{x}^r \J{\xi}^{\rho} V_{g}f(x,\xi) \|\J{\eta}^{|\rho|} |V_{\psi}[V_2(t,\cdot)](z,\eta)\|_{L_{z}^{\infty}} dxdy d\xi d\eta  \notag
\\
&\lesssim \|V_2(t)\|_{M_{0,|\rho|}^{\infty,1}} \int_{\R^{2d}} \left(\int \widetilde{g}(y-x)dy\right)  \J{x}^r \J{\xi}^{\rho} |V_{g}f(x,\xi)| dxd\xi \notag
\\
&\lesssim \|f\|_{W_{r,\rho}^{1,1}}.
\end{align}
Due to (\ref{inq6}) and (\ref{inq7}), we derive
\begin{equation} \label{inq8}
\|K_2(t,s)f\|_{W_{r,\rho}^{1,1}}\lesssim \|f\|_{W_{r,\rho}^{1,1}}.
\end{equation}
Similarly, we have also
\begin{equation} \label{inq9}
\|K_2(t,s)f\|_{W_{r,\rho}^{p,q}}\lesssim \|f\|_{W_{r,\rho}^{p,q}}
\end{equation}
with $(p,q)=(1,\infty), (\infty,1), (\infty, \infty)$.
Hence it turns out from (\ref{inq8}), (\ref{inq9}) and Proposition \ref{Prop4} that (\ref{Kv2}) holds.
%%%%%%%%%%%%%%%%%%%%%%%%%%%%%%%%%%%%%%%%%%%%%%%%%
%%%%%%%%%%%%%%%%%%%%%%%%%%%%%%%%%%%%%%%%%%%%%%%%%
\subsection{Proof of well-posedness}
Let $K(t,s)=K_{\sigma}(t,s)+K_{V_1}(t,s)+K_{V_2}(t,s)$ for $0\leq|s-t|\leq T$.
By the argument in Section 3.2, the operator $K(t,s)$ is uniformly bounded on $W_{r,\rho}^{p,q}$ with respect to $0\leq|s-t|\leq T$ for all $1\leq p, q\leq\infty$ and $r, \rho \in \R$.
Hence there exists a solution $v\in C([0,T]; \CAL{W}_{r,\rho}^{p,q}(\R^d;\C^n))$ to
\begin{equation}\label{IEv}
v(t)=K(t,0)u_0+i\int_{0}^{t}K(t,s)v(s) ds
\end{equation}
by the Picard iteration. Moreover Gronwall's inequality follows that
\begin{equation}\label{vbdd}
\|v(t)\|_{W_{r,\rho}^{p,q}}\leq C_T \|u_0\|_{W_{r,\rho}^{p,q}}
\end{equation}
with some constant $C_T$ depending on $T$. 
We set 
\begin{equation}
U(t)u_0=U_1(t,0)u_0+i\int_{0}^{t}U_1(t,s) v(s) ds
\end{equation}
for the solution $v\in C([0,T]; \CAL{W}_{r,\rho}^{p,q}(\R^d;\C^n))$ to (\ref{IEv}).
Then by (\ref{unitary}), (\ref{DE}) and (\ref{IEv}), $U(t)u_0$ solves (\ref{CP2}) in $ C([0,T]; \CAL{W}_{r,\rho}^{p,q}(\R^d;\C^n))$.
Moreover by (\ref{inq4}) and (\ref{vbdd}), we have 
\begin{equation}
\|U(t)u_0\|_{W_{r,\rho}^{p,q}}\leq C_T \|u_0\|_{W_{r,\rho}^{p,q}}, 
\end{equation} 
which shows uniqueness of solutions to (\ref{CP2}) in $C([0,T]; \CAL{W}_{r,\rho}^{p,q}(\R^d;\C^n))$ and  the solution is continuous on $\CAL{W}_{r,\rho}^{p,q}(\R^d;\C^n)$ with respect to initial data.
%%%%%%%%%%%%%%%%%%%%%%%%%%%%%%%%%%%%%%%%%%%%%%%%%
%%%%%%%%%%%%%%%%%%%%%%%%%%%%%%%%%%%%%%%%%%%%%%%%%
\section{Proof of Theorem \ref{thm2}}
We set
\begin{align*}
U_2(t,s)f(y)&=\iint e^{-i\int_{s}^{t}\left( \sigma(\xi(\tau))+V(\tau,x)\right)d\tau}e^{i(y-x)\xi(t)}
 g(y-x)V_{g}f(x,\xi(s))dxd\xi
\end{align*}
for $f\in \CAL{S}(\R^d; \C^n)$ and $g\in \CAL{S}(\R^d)\setminus \{0\}$ with $\|g\|_{L^2}=1$,
where $\xi(t)=\xi-\int_{0}^{t}\nabla_{x}Q(\tau,x)d\tau$.
It follows from the second order Taylor expansion of $V(t,y)$ at $y=x$ that
\begin{align*}
&\left(i\partial_{t}-\sigma(D_{y})-V(t,y)\right)U_2(t,s)f(y)
\\
&=-\iint e^{-i\int_{s}^{t}\left( \sigma(\xi(\tau))+V(\tau,x)\right)d\tau}e^{i(y-x)\xi(t)} V_{g}f(x,\xi(s))(R_{V}(t,x,y)+R_{\sigma'}(t,x,y,\xi)) dx d\xi, 
\end{align*}
where
\begin{align*}
&R_{V}(t,x,y)=\frac{1}{2}\sum_{|\alpha|=2} \int_{0}^{1}\left(\partial_{x}^{\alpha}V\right)(t,x+\theta(y-x)) (1-\theta)d\theta (y-x)^{\alpha} g(y-x), 
\\
&R_{\sigma'}(t,x,y,\xi) =\int e^{(y-x)\eta} \left( \sigma(\eta+\xi(t))-\sigma(\xi(t)) \right) \widehat{g} (\eta) d\eta.
\end{align*}
Let
\begin{align*}
&K_{\sigma'}(t,s)f(y)=-\iint e^{-i\int_{s}^{t}\left( \sigma(\xi(\tau))+V(\tau,x)\right)d\tau}e^{i(y-x)\xi(t)} V_{g}f(x,\xi(s)) R_{\sigma'}(t,x,y,\xi) dx d\xi,
\\
&K_{V}(t,s)f(y)= -\iint e^{-i\int_{s}^{t}\left( \sigma(\xi(\tau))+V(\tau,x)\right)d\tau}e^{i(y-x)\xi(t)} V_{g}f(x,\xi(s)) R_{V}(t,x,y) dx d\xi.
\end{align*}
Then we have
\begin{equation}\label{DE2}
\left(i\partial_{t}-\sigma(D_{y})-V(t,y)\right)U_0(t,s)f(y)=K_{\sigma'}(t,s)f(y)+K_{V}(t,s)f(y).
\end{equation}
We show that operators $U_2(t,s), K_{\sigma'}(t,s)$ and $K_{V}(t,s)$ are uniformly bounded on $W_{r,0}^{p,q}$ with respect to $0\leq|s-t|\leq T$ for all $1\leq p, q\leq\infty$ and $r \in \R$.
First we prove 
\begin{equation} \label{inqU2}
\|U_2(t,s)f\|_{W_{r,0}^{p,q}}\lesssim \|f\|_{W_{r,0}^{p,q}}.
\end{equation}
Let $\varphi \in \CAL{S}(\R^n)\setminus \{0\}$ be fixed.
Using the equality $(1-\Delta_{y})^{N}e^{iy(\xi(t)-\Xi)}=\J{\xi(t)-\Xi}^{2N}e^{iy(\xi(t)-\Xi)}$ and changing variable : $\xi \mapsto \xi(s)$, we have for any $N\in \N,$
\begin{align*}
&|V_{\varphi}[U_2(t,s)f](X,\Xi)|
\\
&=\left|\int_{\R^{3d}}\overline{\varphi(y-X)}e^{-i\int_{s}^{t}\left( \sigma(\xi(\tau))+V(\tau,x)\right)d\tau}e^{i(y-x)\xi(t)-iy\Xi} g(y-x) V_{g}f(x,\xi(s))dxdyd\xi \right|
\\
&\lesssim \sum_{|\beta_1+\beta_2|\leq 2N} \int_{\R^{3d}} \J{\xi(t)-\Xi}^{-2N} |\partial_{y}^{\beta_1}\varphi(y-X)| |\partial_{y}^{\beta_2}g(y-x)| |V_{g}f(x,\xi(s))| dx dy d\xi
\\
&= \sum_{|\beta_1+\beta_2|\leq 2N} \int_{\R^{3d}} \J{\xi-\Xi+\int_{s}^{t}\nabla_{x}Q(\tau,x)d\tau}^{-2N} |\partial_{y}^{\beta_1}\varphi(y-X)| |\partial_{y}^{\beta_2}g(y-x)| |V_{g}f(x,\xi)| dx dy d\xi.
\end{align*}
We write $\widetilde{\varphi}(x)=\J{x}^{|r|} \partial^{\beta_1}\varphi(x)$ and $\widetilde{g}(x)=\J{x}^{|r|} \partial^{\beta_2}g(x)$. Then by Peetre's inequality we obtain
\begin{align*}
&\|\|\J{X}^r V_{\varphi}[U_2(t,s)f](X,\Xi)\|_{L_{\Xi}^1} \|_{L_{X}^1}
\\
&\lesssim \sum_{|\beta_1+\beta_2|\leq 2N} \left\| \left\| \int_{\R^{3d}} \J{\xi-\Xi+\int_{s}^{t}\nabla_{x}Q(\tau,x)d\tau}^{-2N} |\widetilde{\varphi}(y-X) \widetilde{g}(y-x)| \right. \right.
\\
&\left. \left. \hspace{40mm} \times \J{x}^r |V_{g}f(x,\xi)| dx dy d\xi \right\|_{L_{\Xi}^1}  \right\|_{L_{X}^1}
\\
&\leq \sum_{|\beta_1+\beta_2|\leq 2N} \int_{\R^{3d}} \|\J{\cdot}^{-2N} \|_{L^1} \| \widetilde{\varphi}\|_{L^1} |\widetilde{g}(y-x)| \J{x}^r |V_{g}f(x,\xi)| dx dy d\xi 
\\
& \lesssim \|\|\J{x}^r V_{g}f(x,\xi)\|_{L_{\xi}^1} \|_{L_{x}^1},
\end{align*}
which implies (\ref{inqU2}) with $(p,q)=(1,1)$. Similary we have
\begin{equation} 
\|U_2(t,s)f\|_{W_{r,0}^{p,q}} \lesssim \|f\|_{W_{r,0}^{p,q}}
\end{equation}
with $(p,q)=(1,\infty), (\infty,1), (\infty, \infty)$.
Proposition \ref{Prop4} follows (\ref{inqU2}).
Combining the argument in (\ref{Ksigma}) and (\ref{inqU2}),  we can prove 
\begin{equation}
\|K_{\sigma'}(t,s)f\|_{W_{r,0}^{p,q}} \lesssim \|f\|_{W_{r,0}^{p,q}}.
\end{equation}
Next we show 
\begin{equation} \label{inqKv}
\|K_{V}(t,s)f\|_{W_{r,0}^{p,q}} \lesssim \|f\|_{W_{r,0}^{p,q}}.
\end{equation}
We write $V^{\alpha}(t,x,y)=\int_{0}^{1}\left(\partial_{x}^{\alpha}V\right)(t,x+\theta(y-x)) (1-\theta)d\theta$ and $g^{\alpha}(x)=x^{\alpha} g(x)$ for short.
Then $R_{V}(t,x,y)=\sum_{|\alpha|=2}V^{\alpha}(t,x,y) g^{\alpha}(y-x)$ is valid.
Inverse formula for STFT and integration by parts show that for $\psi\in \CAL{S}(\R^n)$ with $\|\psi\|_{L^2}=1$,
\begin{align} \label{inq15}
&|V_{\varphi}[K_V(t,s)f](X,\Xi)| \notag
\\
&\leq \sum_{|\alpha|=2} \left| \int_{\R^{3d}} \overline{\varphi(y-X)} e^{-i\int_{s}^{t}\left( \sigma(\xi(\tau))+V(\tau,x)\right)d\tau}e^{i(y-x)\xi(t)-iy\Xi} g^{\alpha}(y-x) V^{\alpha}(t,x,y) V_{g}f(x,\xi(s))dxdyd\xi \right| \notag
\\
&=\sum_{|\alpha|=2} \left| \int_{\R^{5d}} \overline{\varphi(y-X)} g^{\alpha}(y-x) \psi(y-z) e^{-i\int_{s}^{t}\left( \sigma(\xi(\tau))+V(\tau,x)\right)d\tau-ix\xi}e^{iy(\xi(t)+\eta-\Xi)} \right. \notag
\\
&\left. \hspace{25mm} \times V_{\psi}[V^{\alpha}(t,x,\cdot)](z,\eta) V_{g}f(x,\xi(s))dxdydzd\xi d\eta \right| \notag
\\
&\lesssim \sum_{|\alpha|=2} \sum_{|\beta_1+\beta_2+\beta_3|\leq 2N} \int_{\R^{5d}} \J{\xi(t)+\eta-\Xi}^{-2N} \left| \partial_{y}^{\beta_1}\overline{\varphi(y-X)} \partial_{y}^{\beta_2}g^{\alpha}(y-x) \partial_{y}^{\beta_3}\psi(y-z) \right| \notag
\\
& \hspace{40mm} \times \left|V_{\psi}[V^{\alpha}(t,x,\cdot)](z,\eta) V_{g}f(x,\xi(s)) \right| dxdydzd\xi d\eta .
\end{align}
Combining the fact $\J{X}^{r}\lesssim \J{y-X}^{|r|} \J{y-x}^{|r|} \J{x}^{r}$ and (\ref{inq5}), we observe
\begin{align} \label{inq16}
&\J{X}^r  |V_{\varphi}[K_2(t,s)f](X,\Xi)| \notag
\\
&\lesssim \sum_{|\alpha|=1} \sum_{|\beta_1+\beta_2+\beta_3|\leq 2N} \int_{\R^{5d}} \J{\xi(t)+\eta-\Xi}^{-2N} \left| \widetilde{\varphi}(y-X) \widetilde{g^{\alpha}}(y-x) \widetilde{\psi}(y-z) \right| \notag
\\
& \hspace{40mm} \times |V_{\psi}[V^{\alpha}(t,x,\cdot)](z,\eta)| \J{x}^r |V_{g}f(x,\xi(s)) | dxdydzd\xi d\eta,
\end{align}
where $\widetilde{\varphi}(x)=\J{x}^{|r|}\partial_{x}^{\beta_1}\varphi(x)$, $\widetilde{g^{\alpha}}(x)=\J{x}^{|r|}\partial_{x}^{\beta_2}g^{\alpha}(x)$ and $\widetilde{\psi}(x)=\partial_{x}^{\beta_3}\psi(x)$.
Translating : $\xi \mapsto \xi(s)$, we have 
\begin{align} \label{inq17}
&\left\| \left\|\int_{\R^{5d}} \J{\xi+\eta-\Xi-\int_{s}^{t}\nabla Q(\tau,x)d\tau}^{-2N} \left| \widetilde{\varphi}(y-X) \widetilde{g^{\alpha}}(y-x) \widetilde{\psi}(y-z) \right| \right.\right. \notag
\\
& \hspace{20mm} \times \left.\left.  |V_{\psi}[V^{\alpha}(t,x,\cdot)](z,\eta)| \J{x}^r |V_{g}f(x,\xi) | dxdydzd\xi d\eta \right\|_{L_{\Xi}^{1}} \right\|_{L_{X}^{1}} \notag
\\
&\lesssim \int_{\R^{4d}} \left\|\J{\cdot}^{-2N} \right\|_{L^1}  \left\| \widetilde{\varphi} \right\|_{L^1} |\widetilde{g^{\alpha}}(y-x)| \J{x}^r |V_{g}f(x,\xi) |   \notag
\\
& \hspace{20mm} \times  V_{\psi}[V^{\alpha}(t,x,\cdot)](z,\eta)\|_{L_{z}^{\infty}} \left( \int |\widetilde{\psi}(y-z)| dz \right) dxdyd\xi d\eta \notag
\\
&\lesssim \int_{\R^{2d}}  \left(\int \widetilde{g^{\alpha}}(y-x) dy \right) \J{x}^r |V_{g}f(x,\xi) | \|\| V_{\psi}[V^{\alpha}(t,x,\cdot)](z,\eta)\|_{L_{z}^{\infty}} \|_{L_{\eta}^{1}}dxdyd\xi  \notag
\\
& \lesssim \sup_{x}\|V^{\alpha}(t,x,\cdot)\|_{M^{\infty,1}} \|f\|_{W_{r,0}^{1,1}} \lesssim \|f\|_{W_{r,0}^{1,1}}.
\end{align}
Due to (\ref{inq16}) and (\ref{inq17}), we derive
\begin{equation} \label{inq18}
\|K_V(t,s)f\|_{W_{r,0}^{1,1}}\lesssim \|f\|_{W_{r,0}^{1,1}}.
\end{equation}
Similarly, we have also
\begin{equation} \label{inq19}
\|K_V(t,s)f\|_{W_{r,0}^{p,q}}\lesssim \|f\|_{W_{r,0}^{p,q}}
\end{equation}
with $(p,q)=(1,\infty),(\infty,1)$ and $(\infty,\infty)$.
Hence it turns out from (\ref{inq18}), (\ref{inq19}) and Proposition \ref{Prop4} that (\ref{inqKv}) holds.
Therefore we complete the proof of Theorem \ref{thm2} by the same argumennt as one in Section 3.3.

%%%%%%%%%%%%%%%%%%%%%%%%%%%%%%%%%%%%%%%%%%%%%%%%%%%%%%%%%%%%%%%%%%%%%%%%%%%%%%%%%%%%%%
%%%%%%%%%%%%%%%%%%%%%%%%%%%%%%%%  $5 参考文献  %%%%%%%%%%%%%%%%%%%%%%%%%%%%%%%%%%%%%%%%%
%%%%%%%%%%%%%%%%%%%%%%%%%%%%%%%%%%%%%%%%%%%%%%%%%%%%%%%%%%%%%%%%%%%%%%%%%%%%%%%%%%%%%%
%%%%%%%%%%%%%%%%%%%%%%%%%%%%%%%%%%%%%%%%%%%%%%%%%%%%%%%%%%%%%%%%%%%%%%%%%%%%%%%%%%%%%%

\address{
Shun Takizawa\\
Department of Mathematics,\\
 Faculty of Science, \\
Tokyo University of Science,\\
Kagurazaka 1-3, Shinjuku-ku, \\
Tokyo 162-8601, Japan}
{1123703@ed.tus.ac.jp}

\end{document}